# Physics-Informed Kernel Function Neural Networks for Solving Partial Differential Equations

Zhuojia Fu[a,b,*], Wenzhi Xu[b], Shuainan Liu[b]

[a] Key Laboratory of Ministry of Education for Coastal Disaster and Protection, Hohai University, Nanjing 210098, China
[b] College of Mechanics and Materials, Hohai University, Nanjing 211100, China
∗Corresponding author, E-mail: paul212063@hhu.edu.cn

**Abstract:** This paper proposed a novel radial basis function neural network (RBFNN) to solve various partial differential equations (PDEs). In the proposed RBF neural networks, the physics-informed kernel functions (PIKFs), which are derived according to the governing equations of the considered PDEs, are used to be the activation functions instead of the traditional RBFs. Similar to the well-known physics-informed neural networks (PINNs), the proposed physics-informed kernel function neural networks (PIKFNNs) also include the physical information of the considered PDEs in the neural network. The difference is that the PINNs put this physical information in the loss function, and the proposed PIKFNNs put the physical information of the considered governing equations in the activation functions. By using the derived physics-informed kernel functions satisfying the considered governing equations of homogeneous, nonhomogeneous, transient PDEs as the activation functions, only the boundary/initial data are required to train the neural network. Finally, the feasibility and accuracy of the proposed PIKFNNs are validated by several benchmark examples referred to high-wavenumber wave propagation problem, infinite domain problem, nonhomogeneous problem, long-time evolution problem, inverse problem, spatial structural derivative diffusion model, and so on.

**Keywords:** Radial basis function neural network; Physics-informed kernel function; Meshless; Activation function

# 1. Introduction

In the past few decades, numerical simulation becomes more and more important in many scientific and engineering fields. As listed in Table 1, the finite element method (FEM) (Moes, Dolbow, & Belytschko, 1999; Tornabene, Fantuzzi, Ubertini, & Viola, 2015) and finite difference method (FDM) (Santoro & Elishakoff, 2006) are two predominant numerical discretization methods for various scientific and engineering applications. However, they still face time-consuming process, laborious mesh generation, and remeshing in the numerical solution of several specific problems, such as wave propagation under infinite domain, large-scale-ratio structures, moving boundary problems, and so on. To overcome these drawbacks, a group of meshless and particle methods have been proposed and successfully applied to the aforementioned problems (J. T. Chen, Chen, Chen, & Chen, 2007; L. Chen, Xu, & Fu, 2022; W. Chen, Lin, & Wang, 2011; Fu, Yang, Zhu, & Xu, 2019; Lin, Chen, Chen, & Jiang, 2013; Reutskiy & Lin, 2018; Sarler, 2009; Sun, Wei, & Chen, 2020; F. Z. Wang & Zheng, 2016). In both mesh-based and meshless methods, the traditional polynomial basis functions are widely used to approximate the solutions of the considered problems. To improve the computational efficiency, both mesh-based and meshless methods employ the basis functions containing the physical and mechanical information of the considered problems to replace the traditional polynomial basis functions, which forms a group of physics-informed kernel function methods, such as boundary element method (Atroshchenko, Hale, Videla, Potapenko, & Bordas, 2017; Shaaban, Anitescu, Atroshchenko, & Rabczuk, 2023), hybrid Trefftz FEMs (Cao, Qin, & Yu, 2013), scaled boundary finite element method (Song & Wolf, 1999), semi-analytical collocation methods (Fu, et al., 2023; Gaspar, 2013; Kolodziej & Grabski, 2018; W. Li, Chen, & Pang, 2016; Zheng & Ma, 2012) and localized semi-analytical collocation methods (Chkadua, Mikhailov, & Natroshvili, 2017; Fu, et al., 2022; Lin, Qiu, & Wang, 2022; F. Wang, Gu, Qu, & Zhang, 2020) and so on.

**Table 1** Bibliographic database search based on the web of science (Date: 2023-05-08).

| Numerical methods | Search phrase in topic field | No. of entries |
|---|---|---|
| FEM | 'finite element(s)' | 610,949 |
| FDM | 'finite difference(s)' | 106,078 |
| FVM | 'finite volume(s)' | 40,996 |
| BEM | 'boundary element(s)' or 'boundary integral(s)' | 34,842 |
| DEM | 'discrete element(s)' | 21,504 |
| Meshless and Particle methods | 'collocation method(s)' or 'meshless' or 'meshfree' or 'element-free' or 'smoothed particle hydrodynamics' or 'material point method' | 30,196 (11,781) |
| Neural Networks | 'neural network(s)' or 'physics-informed neural network(s)' or 'radial basis function neural network(s)' | 686,665 (9,361) |

On the other hand, the neural network (Hinton & Salakhutdinov, 2006) is considered a powerful tool to be widely used in the fields of computer vision, biomedicine, and oil and gas engineering, resulting in a technological revolution in varied science and engineering fields. Recently, deep learning neural networks (Yann LeCun, Bengio, & Hinton, 2015) become a new research hotspot due to its strong learning ability, which can not only discover physical laws but also solve partial differential equations. Physics-informed neural network (PINN) (Mao, Jagtap, & Karniadakis, 2020; Raissi, Perdikaris, & Karniadakis, 2019) is one of the most representative deep learning neural networks for solving PDEs, which is a standard feedforward neural network associated with the redefined loss function including the physical information from PDEs and initial/boundary conditions. Although the PINN has been successfully used in many scientific and engineering fields (Goswami, Anitescu, Chakraborty, & Rabczuk, 2020; Karniadakis, et al., 2021; Lu, Meng, Mao, & Karniadakis, 2021; Raissi, et al., 2019; Robinson, Pawar, Rasheed, & San, 2022; Roy, Bose, Sundararaghavan, & Arroyave, 2023), it still faces several challenging problems such as extra computational cost for solving linear PDEs, low computational efficiency/difficulty to compute for long-time evolution problems, high wavenumber problems, infinite domain problems. To save the computational cost for solving linear PDEs, the two-layer neural network can be used in the PINN, in which the selection of the neuron kernel functions plays an important role in its ability to approximate the PDE solution. Based on the dimensionality independence and simple form of radial basis functions (RBFs) (Karageorghis & Chen, 2022; Tatari & Dehighan, 2009), a kind of neural network that takes radial basis functions as neuron basis functions, named RBF neural networks (Gorbachenko & Zhukov, 2017; J. Y. Li, Luo, Qi, & Huang, 2003) have been successfully used to solve varied PDEs several decades before PINN was proposed. However, similar to PINN, the standard RBFNNs also have the issue of low computational efficiency/difficulty to compute for solving the above-mentioned problems.

In this study, we combine the advantages of PINNs and RBFNNs and use the physics-informed kernel function containing the considered governing equation information as the neuron basis function to build a class of two-layer physics-informed kernel function neural network (PIKFNN). In comparison with the RBFNNs, the PIKFNN employs the physics-informed kernel function as the neuron basis function instead of the traditional RBFs. In comparison with the PINNs, the PIKFNN employs a two-layer neural network instead of the multiple-layer and introduces the neuron basis function containing the considered governing equation information instead of the constraint condition of the governing equation in the loss function. A brief outline of the paper is given as follows. Section 2 introduces the construction of the PIKFs. In Section 3, the constructions of the PIKFNNs for several PDEs are presented. Section 4 investigates the feasibility and accuracy of the PIKFNNs under some benchmark examples. Finally, some conclusions are summarized in Section 5.

## 2. Physics-informed kernel functions

Most physical and engineering problems can be simplified to partial differential equations (PDEs) subjected to boundary conditions. Without loss of generality, the following elliptical PDE problems are considered

$$L_0\mathbf{u}(\mathbf{x}) = \mathbf{f}(\mathbf{x}), \quad \mathbf{x} \in \Omega, \quad (1)$$

$$\mathbf{Bu}(\mathbf{x}) = \mathbf{g}(\mathbf{x}), \quad \mathbf{x} \in \partial\Omega, \quad (2)$$

where $L_0$, $\mathbf{B}$ represent differential governing equation operator and boundary condition equation operator, $\mathbf{u}(\mathbf{x})$ is the physical quantities at node $\mathbf{x}$ of the considered domain $\Omega$, $\mathbf{f}(\mathbf{x})$, $\mathbf{g}(\mathbf{x})$ denote the prescribed source function inside $\Omega$ and the prescribed function on the boundary $\partial\Omega$.

2.1 Homogeneous PDEs

When the source function $\mathbf{f}(\mathbf{x}) = \mathbf{0}$, Eq. (1) becomes the homogeneous governing equation

$$L_0\mathbf{u}(\mathbf{x}) = \mathbf{0}, \quad \mathbf{x} \in \Omega, \quad (3)$$

for some specific differential governing equation operators $L_0$, the physics-informed kernel functions, such as the fundamental solutions, harmonic functions, radial Trefftz functions and T-complete functions, are listed in Tables 2-4, which satisfy as a prior Eq. (3). In the tables, $\Delta$, $\nabla$ stand for Laplace operator and gradient operator, $D$ denotes the diffusion coefficient, $k$ represents wave number, $\mathrm{i} = \sqrt{-1}$, $J_n$ and $Y_n$ denote $n$-th Bessel function of the first and second kinds, $I_n$ and $K_n$ denote $n$-th modified Bessel function of the first and second kinds, $H_n^{(1)}$ denote $n$-th Hankel function of the first kind, $\mathbf{v}$ and $\mathbf{r} = (r_1, \cdots, r_{\dim})$ denote the vectors of velocity and distance, in which dim is the problem dimension, $r$ is the Euclidean distance between two nodes, $\mu = \sqrt{\left(\frac{|\mathbf{v}|}{2D}\right)^2 + \frac{k}{D}}$, $c$ represents the shape parameter, whose proper value is related to the characteristic length of the computing domain. $P_v^m$ represents associated Legendre polynomials. In 2D problems, $\rho = \sqrt{x_1^2 + x_2^2}$, $\theta = \arctan(x_2/x_1)$; in 3D problems, $\rho = \sqrt{x_1^2 + x_2^2 + x_3^2}$, $\varphi = \arccos(x_3/\sqrt{x_1^2 + x_2^2})$, $\theta = \arctan(x_2/x_1)$.

**Table 2** Fundamental solutions $G_F$ of common-used differential equation operators.

| $\Re$ | 2D | 3D | 4D |
| --- | --- | --- | --- |

| $\Re$ | | | |
|---|---|---|---|
| $\Delta$ | $-\dfrac{\ln(r)}{2\pi}$ | $\dfrac{1}{4\pi r}$ | $\dfrac{1}{4\pi^2 r^2}$ |
| $\Delta + k^2$ | $\dfrac{iH_0^{(1)}(kr)}{4}$ | $\dfrac{e^{-ikr}}{4\pi r}$ | / |
| $\Delta - k^2$ | $\dfrac{K_0(kr)}{2\pi}$ | $\dfrac{e^{-kr}}{4\pi r}$ | / |
| $D\Delta + \mathbf{v}\bullet\nabla - k^2$ | $\dfrac{K_0(\mu r)e^{-\frac{\mathbf{v}\bullet\mathbf{r}}{2D}}}{2\pi}$ | $\dfrac{e^{-\mu r - \frac{\mathbf{v}\bullet\mathbf{r}}{2D}}}{4\pi r}$ | / |
| $\Delta^2$ | $\dfrac{r^2\ln(r)-r^2}{8\pi}$ | $\dfrac{r}{8\pi}$ | / |

**Table 3** Harmonic functions $G_H$ and radial Trefftz functions $G_{RT}$ of common-used differential equation operators.

| $\Re$ | 2D | 3D |
|---|---|---|
| $\Delta$ | $e^{-c(r_1^2-r_2^2)}\cos(2cr_1r_2)$ | $e^{-c(r_1^2-r_2^2)}\cos(2cr_1r_2)$ $+e^{-c(r_2^2-r_3^2)}\cos(2cr_2r_3)$ $+e^{-c(r_3^2-r_1^2)}\cos(2cr_3r_1)$ |
| $\Delta + k^2$ | $\dfrac{1}{2\pi}J_0(kr)$ | $\dfrac{\sin kr}{4\pi r}$ |
| $\Delta - k^2$ | $\dfrac{1}{2\pi}I_0(kr)$ | $\dfrac{\sinh(kr)}{4\pi r}$ |
| $D\Delta + \mathbf{v}\bullet\nabla - k^2$ | $\dfrac{1}{2\pi}I_0(\mu r)e^{-\frac{\mathbf{v}\bullet\mathbf{r}}{2D}}$ | $\dfrac{\sinh(\mu r)}{4\pi r}e^{-\frac{\mathbf{v}\bullet\mathbf{r}}{2D}}$ |
| $\Delta^2$ | $r^2 e^{-c(r_1^2-r_2^2)}\cos(2cr_1r_2)$ | $r^2\Big(e^{-c(r_1^2-r_2^2)}\cos(2cr_1r_2)$ $+e^{-c(r_2^2-r_3^2)}\cos(2cr_2r_3)$ $+e^{-c(r_3^2-r_1^2)}\cos(2cr_3r_1)\Big)$ |

**Table 4** T-complete functions $G_T$ of common-used differential equation operators.

| $\Re$ | 2D | 3D |
|---|---|---|
| $\Delta$ | $1$ | $P_v^0(\cos\varphi)$ |
| | $\rho^m\cos(m\theta)$ | $\rho^m P_v^m(\cos\varphi)\cos(m\theta)$ |
| | $\rho^m\sin(m\theta)$ | $\rho^m P_v^m(\cos\varphi)\sin(m\theta)$ |
| $\Delta + k^2$ | $J_0(k\rho)$ | $J_0(k\rho)P_v^0(\cos\varphi)$ |
| | $J_m(k\rho)\cos(m\theta)$ | $J_m(k\rho)P_v^m(\cos\varphi)\cos(m\theta)$ |
| | $J_m(k\rho)\sin(m\theta)$ | $J_m(k\rho)P_v^m(\cos\varphi)\sin(m\theta)$ |
|  | $I_0(k\rho)$ | $I_0(k\rho)P_v^0(\cos\varphi)$ |
| $\Delta - k^2$ | $I_m(k\rho)\cos(m\theta)$ | $I_m(k\rho)P_v^m(\cos\varphi)\cos(m\theta)$ |
|  | $I_m(k\rho)\sin(m\theta)$ | $I_m(k\rho)P_v^m(\cos\varphi)\sin(m\theta)$ |
|  | $\rho^2$ | $\rho^2 P_v^0(\cos\varphi)$ |
| $\Delta^2$ | $\rho^{m+2}\cos(m\theta)$ | $\rho^{m+2}P_v^m(\cos\varphi)\cos(m\theta)$ |
|  | $\rho^{m+2}\sin(m\theta)$ | $\rho^{m+2}P_v^m(\cos\varphi)\sin(m\theta)$ |

2.2 Nonhomogeneous PDEs

When the specific non-zero source function $f(\mathbf{x})$ is sufficiently smooth, the physics-informed kernel function (PIKF) satisfying Eq. (1) can be composed of two parts, the PIKF satisfying homogeneous Eq. (3) and the PIKFs related to $f(\mathbf{x})$ by satisfying the following high-order homogeneous equation

$$L_{NL}\ldots L_2 L_1(f(\mathbf{x}))=0, \quad \mathbf{x}\in\Omega \tag{4}$$

where $L_1, L_2, \ldots L_{NL}$ represent the same or different differential operators. If all the differential operators $L_0, L_1, L_2, \ldots L_{NL}$ are different, the corresponding PIKFs can be found from Tables 2-4 listed in subsection 2.1. If some of the differential operators $L_0, L_1, L_2, \ldots L_{NL}$ are the same, then the PIKFs for the related high-order differential operator $\Re^n$ need to be employed. Tables 5-7 list these high-order PIKFs, such as the fundamental solutions, harmonic functions, radial Trefftz functions and T-complete functions. In the tables, $A_0=1, A_n=\dfrac{A_{n-1}}{2nk^2}$, $B_0=0$, $B_{n+1}=\dfrac{1}{4(n+1)^2}\left(\dfrac{C_n}{n+1}+B_n\right)$, $C_0=1$, $C_{n+1}=\dfrac{C_n}{4(n+1)^2}$, $D_0=1$, $D_{n+1}=k\rho D_n$.

**Table 5** Fundamental solutions $G_F^n$ of common-used high-order differential equation operators.

| $\Re^n$ | 2D | 3D |
|---|---|---|
| $\Delta^{n+1}$ | $\dfrac{r^{2n}}{2\pi}(C_n \ln r - B_n)$ | $\dfrac{1}{(2n)!}\dfrac{r^{2n-1}}{4\pi}$ |
| $(\Delta+k^2)^{n+1}$ | $A_n(kr)^{n+1-\dim/2}\mathrm{i}H^{(1)}_{n-1+\dim/2}(kr)$ | |
| $(\Delta-k^2)^{n+1}$ | $A_n(kr)^{n+1-\dim/2}K_{n-1+\dim/2}(kr)$ | |
| $\left(D\Delta + \mathbf{v}\bullet\nabla - k^2\right)^{n+1}$ | | $A_n(\mu r)^{n+1-\dim/2} K_{n-1+\dim/2}(\mu r) e^{-\frac{\mathbf{v}\bullet\mathbf{r}}{2D}}$ |

**Table 6** Harmonic functions $G_H^n$ and radial Trefftz functions $G_{RT}^n$ of common-used high-order differential equation operators.

| $\Re^n$ | 2D | 3D |
|---|---|---|
| $\Delta^{n+1}$ | $r^{2n} e^{-c(r_1^2 - r_2^2)} \cos(2cr_1 r_2)$ | $r^{2n}\left(e^{-c(r_1^2-r_2^2)}\cos(2cr_1r_2) + e^{-c(r_2^2-r_3^2)}\cos(2cr_2r_3) + e^{-c(r_3^2-r_1^2)}\cos(2cr_3r_1)\right)$ |
| $\left(\Delta+k^2\right)^{n+1}$ | $A_n(kr)^{n+1-\dim/2} J_{n-1+\dim/2}(kr)$ | |
| $\left(\Delta-k^2\right)^{n+1}$ | $A_n(kr)^{n+1-\dim/2} I_{n-1+\dim/2}(kr)$ | |
| $\left(D\Delta+\mathbf{v}\bullet\nabla-k^2\right)^{n+1}$ | $A_n(\mu r)^{n+1-\dim/2} I_{n-1+\dim/2}(\mu r) e^{-\frac{\mathbf{v}\bullet\mathbf{r}}{2D}}$ | |

**Table 7** T-complete functions $G_{RT}^n$ of common-used high-order differential equation operators.

| $\Re^n$ | 2D | 3D |
|---|---|---|
| $\Delta^{n+1}$ | $\rho^{2n}$ | $\rho^{2n} P_\nu^0(\cos\varphi)$ |
| | $\rho^{m+2n}\cos(m\theta)$ | $\rho^{m+2n} P_\nu^m(\cos\varphi)\cos(m\theta)$ |
| | $\rho^{m+2n}\sin(m\theta)$ | $\rho^{m+2n} P_\nu^m(\cos\varphi)\sin(m\theta)$ |
| $\left(\Delta+k^2\right)^{n+1}$ | $D_n J_n(k\rho)$ | $D_n J_n(k\rho) P_\nu^0(\cos\varphi)$ |
| | $D_n J_{m+n}(k\rho)\cos(m\theta)$ | $D_n J_{m+n}(k\rho) P_\nu^m(\cos\varphi)\cos(m\theta)$ |
| | $D_n J_{m+n}(k\rho)\sin(m\theta)$ | $D_n J_{m+n}(k\rho) P_\nu^m(\cos\varphi)\sin(m\theta)$ |
| $\left(\Delta-k^2\right)^{n+1}$ | $D_n I_n(k\rho)$ | $D_n I_n(k\rho) P_\nu^0(\cos\varphi)$ |
| | $D_n I_{m+n}(k\rho)\cos(m\theta)$ | $D_n I_{m+n}(k\rho) P_\nu^m(\cos\varphi)\cos(m\theta)$ |
| | $D_n I_{m+n}(k\rho)\sin(m\theta)$ | $D_n I_{m+n}(k\rho) P_\nu^m(\cos\varphi)\sin(m\theta)$ |

2.3 Transient PDEs

The differential equation operators discussed in the last two subsections are time-independent operators. From this subsection, we will focus on the time dependent differential operators. Consider the following transient PDE problems,

$$\left(\frac{\partial^m}{\partial t^m} - L_0\right) u(\mathbf{x},t) = f(\mathbf{x},t), \ \mathbf{x}\in\Omega, \ t\in(0,T], \tag{5}$$

$$Bu(\mathbf{x},t) = g(\mathbf{x},t), \ \mathbf{x}\in\partial\Omega, \ t\in(0,T], \tag{6}$$

$$Iu(\mathbf{x},0) = g_I(\mathbf{x}), \ \mathbf{x}\in\Omega, \tag{7}$$

where $I=\left[1,\cdots,\frac{\partial^{m-1}}{\partial t^{m-1}}\right]^T$ denotes the initial condition equation operator, $f(\mathbf{x},t)$ stands for the prescribed time-dependent/independent source function, $g(\mathbf{x},t)$ and $g_I(\mathbf{x})$ are the prescribed functions. The physics-informed kernel function (PIKF) satisfying Eq. (5) can be also composed of two parts, the PIKF satisfying homogeneous transient equation ($f(\mathbf{x},t)=0$) and the PIKFs related to $f(\mathbf{x},t)$ by satisfying the following high-order homogeneous equation

$$\bar{L}_{NL}\ldots\bar{L}_2\bar{L}_1(f(\mathbf{x},t))=0, \mathbf{x}\in\Omega, t\in(0,T], \tag{8}$$

where $\bar{L}_1,\bar{L}_2,\cdots,\bar{L}_{NL}$ represent the same or different time-dependent/independent differential operators. The time-dependent PIKFs, such as the fundamental solutions and radial Trefftz functions, are listed in Tables 8 and 9, which satisfy as a prior homogeneous transient equation.

**Table 8** Fundamental solutions $G_F$ of time-dependent differential equation operators.

| $\partial^m/\partial t^m - \Re$ | 2D | 3D |
|---|---|---|
| $\frac{\partial}{\partial t}-k\Delta$ | $\frac{\Theta(t-\tau)e^{r^2/4k(t-\tau)}}{4\pi k(t-\tau)}$ | $\frac{\Theta(t-\tau)e^{r^2/4k(t-\tau)}}{(4\pi k(t-\tau))^{3/2}}$ |
| $\frac{\partial^2}{\partial t^2}-c_1\Delta$ | $\frac{\Theta(t-r/c_1)}{2\pi c_1\sqrt{c_1^2 t^2 - r^2}}$ | $\frac{\Theta(t-r/c_1)}{4\pi r}$ |

**Table 9** Radial Trefftz functions $G_{RT}$ of time-dependent differential equation operators.

| $\partial^m/\partial t^m - \Re$ | 2D | 3D |
|---|---|---|
| $\frac{\partial}{\partial t}-k\Delta$ | $e^{-k(t-\tau)}J_0(r)$ | $e^{-k(t-\tau)}\frac{\sin(r)}{r}$ |
| $\frac{\partial^2}{\partial t^2}-c_1\Delta$ | $\cos(c_1(t-\tau))J_0(r)$ $+\sin(c_1(t-\tau))J_0(r)$ | $\frac{\cos(c_1(t-\tau))\sin(r)}{r}$ $+\frac{\sin(c_1(t-\tau))\sin(r)}{c_1 r}$ |

2.4 Spatial-temporal structural derivative PDEs

The classical differential equation operators discussed in the last two sections are standard partial derivative operators. However, some recent studies have found that phenomena such as anomalous diffusion and acoustic attenuation based on power-law frequency dependence are difficult to be described by the above classical differential equations. In this section, the spatial-temporal structural derivative operator (F. Wang, Cai, Zheng, & Wang, 2020) is constructed to describe these anomalous physical and mechanical phenomena. The temporal and spatial structural derivative can be defined as

$$\frac{du}{d_s t} = \lim_{t \to t'} \frac{u(t) - u(t')}{\vartheta(t) - \vartheta(t')}, \tag{9}$$

$$\frac{du}{d_s x} = \lim_{x \to x'} \frac{u(x) - u(x')}{F(x) - F(x')}, \tag{10}$$

where $d_s$ stands for the structural derivative, $\vartheta(t)$ and $F(x)$ are arbitrary structural functions. When $\vartheta(t) = t^\alpha$ and $F(\mathbf{x}) = \mathbf{x}^\beta$ are power functions, the structural derivatives are degenerated to Hausdorff derivative ($0 < \alpha, \beta < 2, \alpha \neq 1, \beta \neq 1$) or standard partial derivative ($\alpha = \beta = 1$). Here consider the following 2D spatial-temporal structural derivative PDE,

$$\frac{du(\mathbf{x},t)}{d_s t} = D\left(\frac{d}{d_s x_1}\left(\frac{du(x_1,t)}{d_s x_1}\right) + \frac{d}{d_s x_2}\left(\frac{du(x_2,t)}{d_s x_2}\right)\right), \quad \mathbf{x} \in \Omega, 0 < t \leq T, \tag{11}$$

subjected to the boundary conditions (6) and initial conditions (7). The corresponding fundamental solutions can be represented in the form of

$$G_F(\mathbf{x},\mathbf{s},t,\tau) = \frac{\Theta(t-\tau)e^{\|F(\mathbf{x})-F(\mathbf{s})\|^2/4D(\vartheta(t)-\vartheta(\tau))}}{(4\pi D(\vartheta(t) - \vartheta(\tau)))}, \tag{12}$$

where $\mathbf{s}$, $\tau$ stand for the spatial and temporal reference nodes, $\Theta$ represents the Heaviside function. When $\vartheta(t) = t$ and $F(\mathbf{x}) = \mathbf{x}$, Eq. (11) goes back to the classical model describing Fickian diffusion process with $\langle x^2(t) \rangle \propto t$, the fundamental solutions (12) can be simplified as

$$G_F(\mathbf{x},\mathbf{s},t,\tau) = \frac{\Theta(t-\tau)e^{\|\mathbf{x}-\mathbf{s}\|^2/4D(t-\tau)}}{(4\pi D(t-\tau))}, \tag{13}$$

when $\vartheta(t)=t^\alpha$ and $F(\mathbf{x})=\mathbf{x}$, Eq. (11) becomes the Hausdorff derivative diffusion model with the mean square displacement proportional to a function of time $t$ $\langle x^2(t)\rangle \propto t^\alpha$, the fundamental solutions (12) can be simplified as

$$G_F(\mathbf{x},\mathbf{s},t,\tau) = \frac{\Theta(t-\tau)e^{\|\mathbf{x}-\mathbf{s}\|^2/4D(t^\alpha-\tau^\alpha)}}{(4\pi D(t^\alpha-\tau^\alpha))}, \quad (14)$$

when $\vartheta(t)=t$ and $F(\mathbf{x})=\mathbf{x}^\beta$, Eq. (11) becomes the Hausdorff derivative diffusion model with the mean square displacement proportional to a function of time t $\langle x^2(t)\rangle \propto t^{(3-\beta)/2\beta}$, the fundamental solutions (12) can be simplified as

$$G_F(\mathbf{x},\mathbf{s},t,\tau) = \frac{\Theta(t-\tau)e^{\|\mathbf{x}^\beta-\mathbf{s}^\beta\|^2/4D(\vartheta(t)-\vartheta(\tau))}}{(4\pi D(\vartheta(t)-\vartheta(\tau)))}, \quad (15)$$

## 3. Physics-informed kernel function neural networks

To set the stage for introducing the proposed PIKFNNs, without loss of generality, we first introduce the PINNs. The raise of PINNs is aimed mainly at solving PDEs, the traditional PINN is a standard feedforward neural network, combined with back propagation algorithm (Y. Lecun, Bottou, Bengio, & Haffner, 1998) and automatic differentiation technique (Baydin, Pearlmutter, Radul, & Siskind, 2018). Take the homogeneous elliptical PDE problems (2)-(3) as an example, a standard PINN is shown in Fig. 1, where $\sigma$ is the nonlinear activation function, $\theta$ is the network parameter to be optimized, contains the weight $w$ and bias $b$. The biggest difference between the PINN and the traditional data-driven neural network is the loss function, in the PINNs, the loss function $Loss_{PINN}$ (16) is defined as the mean-square error (MSE) of the residual errors of the boundary conditions and the governing equation,

$$Loss_{PINN} = \frac{1}{N_B}\sum_{i=1}^{N_B}(\mathbf{B}\mathbf{u}_{PINN}(\mathbf{x}_i)-\mathbf{g}(\mathbf{x}_i))^2 + \frac{1}{N_I}\sum_{i=1}^{N_I}(L_0\mathbf{u}_{PINN}(\mathbf{x}_i)-\mathbf{f}(\mathbf{x}_i))^2, \quad (16)$$

where $\mathbf{u}_{PINN}$ is the approximation solution of the PINN, $N_I$ denotes the interior training node number, $N_B$ denotes the boundary training node number, and the total training data number is $N = N_I + N_B$. In the loss function $Loss_{PINN}$, the first term

$\frac{1}{N_B}\sum_{i=1}^{N_B}\left(L_0\mathbf{u}_{PINN}(\mathbf{x}_i)-\mathbf{g}(\mathbf{x}_i)\right)^2$ denotes the data fitting, which is usually used in data-driven neural networks; And the last term $\frac{1}{N_I}\sum_{i=1}^{N_I}\left(L_0\mathbf{u}_{PINN}(\mathbf{x}_i)-\mathbf{f}(\mathbf{x}_i)\right)^2$ denotes the physics regularization, which is the biggest innovation of the PINNs. Based on the chain rule, the automatic differentiation technology is introduced, the differential term of the PINN solution is approximated by automatic differential technique. Then some commonly used optimization algorithms including Adam algorithm and L-BFGS algorithm are employed to minimize the loss function $Loss_{PINN}$, this process is called 'training'. After the PINN is trained well, the trained PINN will output the predicted solution $\mathbf{u}_{PINN}(\mathbf{x})$ at the input $\mathbf{x}$.

Up to now, the PINNs have been successfully used in many scientific and engineering fields (Hu, Guo, Long, Li, & Xu, 2022; Kashefi & Mukerji, 2022; Xu, Han, Cheng, Cheng, & Ge, 2022; Zhao, Stuebner, Lua, Phan, & Yan, 2022), however, they still face the challenges of extra computational cost for solving linear PDEs, low computational efficiency/difficulty to computing for long-time evolution problems, high wavenumber problems and infinite domain problems.

Introducing structured prior information to learning machines can quickly make them converge to the right solution and generalize well when fewer training data are available. Inspired by this point, we build a class of two-layer physics-informed kernel function neural networks (PIKFNNs), which use the physics-informed kernel function (PIKF) containing the prior physical information of considered governing equation as the neuron basis function instead of the traditional nonlinear activation function. Fig. 2 plots the structure of the PIKFNN for solving the homogeneous elliptical PDE problems (2)-(3), where $\mathbf{s}_i$ donates the RBF center, also called source nodes, $p_i$ donates the network parameter, $\varphi$ donates the PIKF. One of the most noticeable structural differences between the PIKFNN and PINN lies in the hidden layer, the PIKFNN uses the single PIKF layer as the hidden layer. Also compared with the PINN, due to the

PIKFNN using the PIKF which naturally contains the governing equation information, the physics regularization naturally exists in the PIKF layer, thus the loss function of the PIKFNN can be written as

$$Loss = \frac{1}{N_B}\sum_{i=1}^{N_B}\left(L_0 \mathbf{u}_{PIKFNN}(\mathbf{x}_i) - \mathbf{g}(\mathbf{x}_i)\right)^2, \quad (17)$$

from which, we can find that the PIKFNN no longer needs the interior node, only needs the boundary training nodes, that is $N = N_B$.

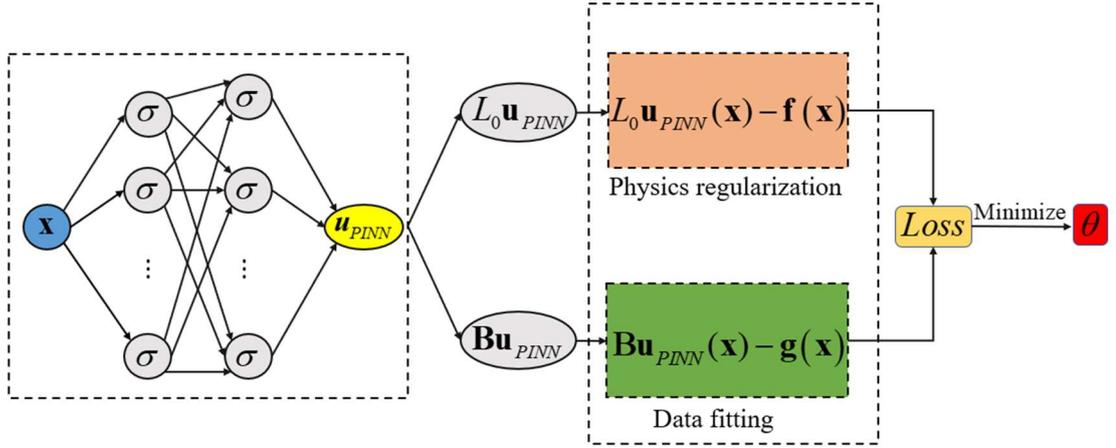

**Fig. 1.** Structure of the PINN for solving the homogeneous PDE problems (2)-(3).

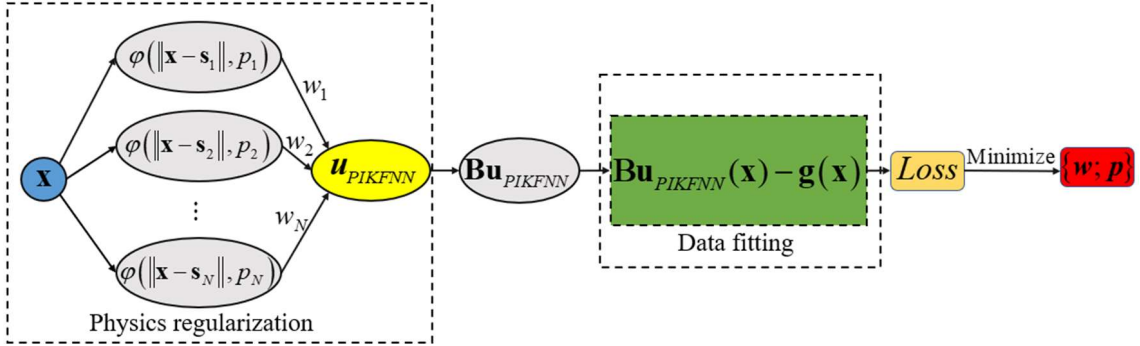

**Fig. 2.** Structure of the PIKFNN for solving the homogeneous PDE problems (2)-(3).

The proposed PIKFNN has the advantages of no interior nodes, simple form, fewer parameters, ease of training, and priori physics regularization. Next, we will introduce the PIKFNN for solving nonhomogeneous PDEs and transient problems.

3.1 PIKFNN for solving nonhomogeneous PDE problems

The PIKFNN for solving the homogeneous problems (2)-(3) has been illustrated above. Here, we will elaborate on the PIKFNN for solving nonhomogeneous problems. To solve the nonhomogeneous problems (1)-(2), some neurons should be added to the PIKF layer as shown in Fig. 3, where $\varphi^{L_0}$ denotes the PIKF corresponding to the governing equation operator $L_0$, and $\varphi^{L_1}$ is the PIKF corresponding to the partial differential operator $L_1$ which satisfies $L_1(f(x))=0, \mathbf{x}\in\Omega$. That is for nonhomogeneous problems the neurons in the PIKF layer are twice the homogeneous problems. Similarly, if you choose a series of $L_{NL}\ldots L_2 L_1$ such that $L_{NL}\ldots L_2 L_1(f(\mathbf{x}))=0, \mathbf{x}\in\Omega$, the additional PIKFs neurons ($\varphi^{L_1},\cdots\varphi^{L_{NL}}$) should be $(NL+1)$ time of the nonhomogeneous problems, where $\varphi^{L_i}, i=1,2,\cdots,NL$ denotes the PIKF corresponding to the differential operator $L_i$.

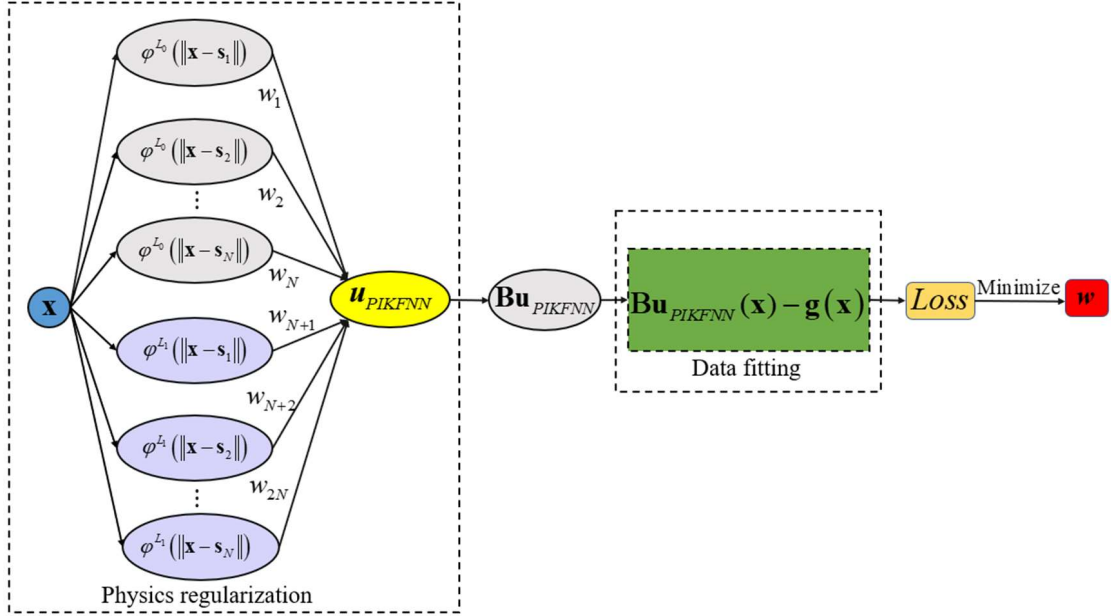

**Fig. 3.** Structure of the PIKFNN for solving the nonhomogeneous PDE problems (1)-(2).

3.2 PIKFNN for solving transient PDE problems

In this subsection, the PIKFNN for solving the transient PDE (5)-(7) is illustrated. As shown in Fig. 4, $\varphi^{\bar{L}_0}$ denotes the PIKF corresponding to the transient governing equation operator, and $\varphi^{\bar{L}_1}$ is the PIKF corresponding to the time-dependent partial differential operator $\bar{L}_1$ which satisfies $\bar{L}_1(f(\mathbf{x},t))=0$, the nodes pair $(\mathbf{s},\tau)$ donates the RBF center in the spatial-temporal domain. The boundary condition and initial condition are used for the data fitting. Corresponding the loss function can be written as

$$Loss = \frac{1}{N_{Bou}}\sum_{i=1}^{N_{Bou}}\left(\mathbf{B}\mathbf{u}_{PINN}(\mathbf{x}_i)-\mathbf{g}(\mathbf{x}_i,t_i)\right)^2 + \frac{1}{N_{Ini}}\sum_{i=1}^{N_{Ini}}\left(\mathbf{I}\mathbf{u}_{PINN}(\mathbf{x}_i)-\mathbf{g}_I(\mathbf{x}_i,t_i)\right)^2, \quad (18)$$

where $N_{Bou}$ donates the boundary condition data number and $N_{Ini}$ donates the initial condition data number, in the PIKFNN for solving transient PDE problems, the total training data number is $N = N_{Ini} + N_{Bou}$.

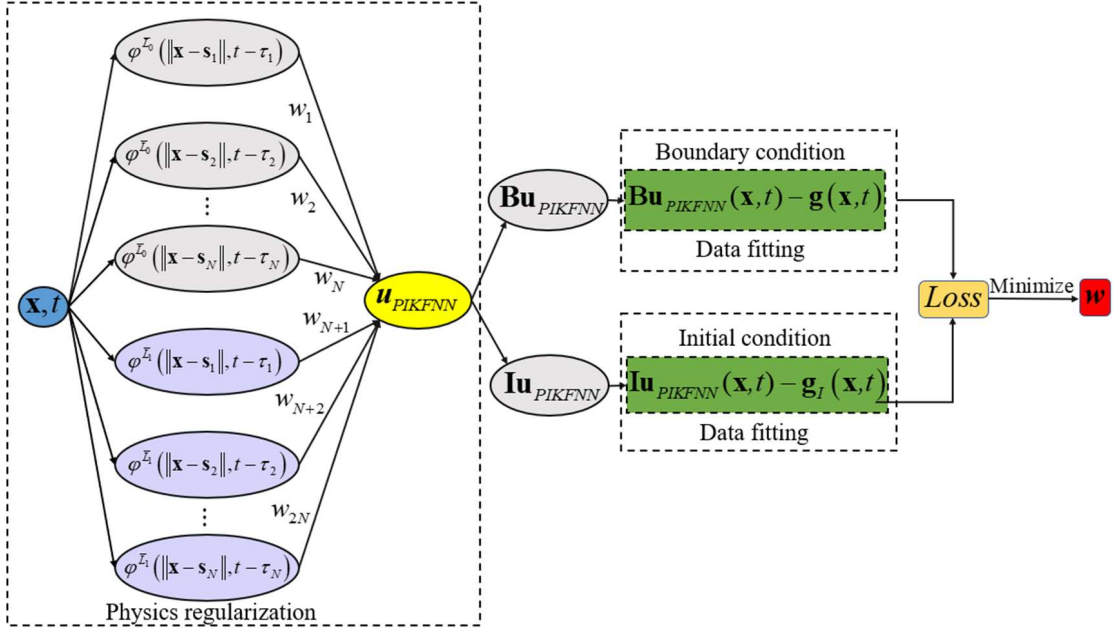

**Fig. 4.** Structure of the PIKFNN for solving the transient PDE problems (5)-(7).

3.3 Training the PIKFNN

The aim of training the PIKFNN is to minimize the loss function,

$$w = \arg\min_{w} Loss \quad (19)$$

In this study, the Adam algorithm and the L-M algorithm are used to minimize the loss function, and the following stopping criterion for the L-M algorithm is used,

*If max |w<sup>i</sup> - w<sup>i-1</sup>| < tol or |Loss<sup>i</sup>(w<sup>i</sup>) - Loss<sup>i-1</sup>(w<sup>i-1</sup>) | < tol,*

***Then** Stop.*

Where the subscript '*i*' denotes the i-*th* iteration, and the *tol* is the tolerance, the tolerance *tol* is a threshold which, if crossed, stops the iterations of the L-M algorithm.

**4. Numerical results and discussions**

To verify the feasibility and accuracy of the proposed PIKFNNs, this section investigates the wave propagation, infinity-Laplace equation, nonhomogeneous modified Helmholtz equation, long-term transient heat conduction equation, spatial structural derivative diffusion, High dimension Laplace equation, potential-based inverse electromyography, elastic problem and so on. To validate its performance, the relative error and the $L_2$ relative error are defined as follows:

$$rerr = (u_{PIKFNN} - u_{Ana})/u_{Ana}, \qquad (20)$$

$$L_2 = \sqrt{\sum_{i=1}^{Nt}(u_{PIKFNN} - u_{Ana})^2 / \sum_{i=1}^{Nt} u_{Ana}^2}, \qquad (21)$$

where $u_{PIKFNN}$ denotes the PIKFNN solution, $u_{Ana}$ represents the analytical solution or measure data, and $Nt$ is the number of test nodes. All the programs are executed by MATLAB software on a desktop with AMD Ryzen 5 5600X 6-Core Processor 3.60 GHz and 16GB RAM.

**Example 1.** Wave propagation in a unit square domain.

Consider the following Helmholtz equation with Dirichlet boundary condition,

$$\left(\Delta + \sqrt{200}^2\right)u(\mathbf{x}) = 0, \quad \mathbf{x} = (x_1, x_2) \in \Omega, \qquad (22)$$

$$u(\mathbf{x}) = \sin(10x_1 + 10x_2), \quad \mathbf{x} = (x_1, x_2) \in \partial\Omega. \qquad (23)$$

The number of boundary training nodes is $N = N_B = 80$, these nodes are uniformly placed on the $[-1,1] \times [-1,1]$ square, the boundary conditions on these nodes are set as the training data for the PIKFNN, in the PIKFNN, the boundary condition values are set as the targets, and their coordinates are set as the labels. The source nodes are uniformly placed on the $[-0.5, 1.5] \times [-0.5, 1.5]$ square, which total number is equal to the boundary training nodes. Due to the problem (22)-(23) is homogeneous, the number of PIKFs was chosen to be $N=80$. The PIKFs corresponding to the partial differential operator $\left(\Delta+\sqrt{200}^2\right)$ is $\varphi^{L_0}=\dfrac{1}{2\pi}Y_0(\sqrt{200}r)$. To minimize the loss function, the Adam algorithm is employed in this example. Unless otherwise stated, in this study, the initialized weights $w$ are set as a series of random numbers which are evenly distributed over the interval from $-1$ to $1$.

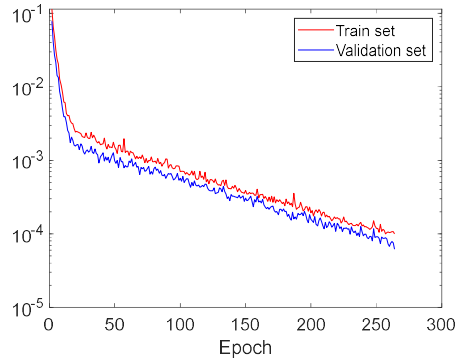

**Fig. 5.** Loss values of the training set and validation set varied with the training epoch in example 1.

First, 80% of the boundary data are set as the training set, and the rest 20% are set as the validation set. From Fig. 5, we can find that the trend of the loss value of the validation set is very similar to the training set, which demonstrates that the validation set is not necessary for the PIKFNNs. Thus, the validation set is no longer presented in the following study.

Next, the convergence of the present PIKFNN is investigated. Table 10 gives the $L_2$ relative errors with different loss values $Loss(w)$, it can be found from Table 10 that with the decreasing Loss value $Loss(w)$, the PIKFNN gives a more accurate predicted solution, correspondingly the CPU time of the training process will increase. This case shows that the proposed PIKFNN is feasible to solve PDEs, and its predicting accuracy depends on the loss value, the smaller the loss value, the higher the predicting accuracy of the PIKFNN.

**Table 10** $L_2$ relative errors with different loss values in Example 1.

| $Loss(\mathbf{w})$ | 1.0E-2 | 1.0E-3 | 1.0E-4 | 1.0E-5 |
|---|---|---|---|---|
| $L_2$ | 7.95E-2 | 2.41E-2 | 2.10E-3 | 1.34E-3 |
| CPU time (s) | 1.10 | 1.58 | 13.2 | 83.8 |

**Example 2.** Wave propagation in a unit circle domain with high wavenumber.

Consider the following Helmholtz equation with Dirichlet boundary condition,

$$(\Delta + k^2)u(\mathbf{x}) = 0 \quad \mathbf{x} = (x_1, x_2) \in \Omega, \tag{24}$$

$$u(\mathbf{x}) = \sin(kx_1) + \cos(kx_2) \quad \mathbf{x} = (x_1, x_2) \in \partial\Omega. \tag{25}$$

The computational domain is a unit circle domain centered at the origin. The source nodes are uniformly placed on a circle of radial 3 centered at the origin. The PIKF corresponding to the partial differential operator $(\Delta + k^2)$ is $\varphi^{L_0} = \frac{1}{2\pi}Y_0(kr)$. To minimize the loss function, the L-M algorithm is employed in this example.

**Table 11** $L_2$ relative errors with different parameter settings in Example 2.

| $k$ | tol | $N_B$ | $L_2$ |
|---|---|---|---|
| 100 | 1E-4 | 300 | 3.27E-7 |
| 100 | 1E-4 | 400 | 1.59E-7 |
| 150 | 1E-4 | 300 | 4.96E-2 |
| 150 | 1E-4 | 400 | 8.96E-6 |
| 150 | 1E-6 | 300 | 2.95E-2 |

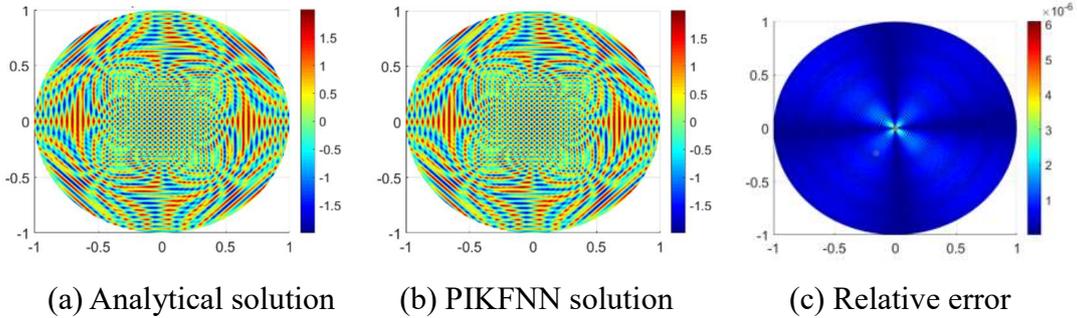

(a) Analytical solution    (b) PIKFNN solution    (c) Relative error

**Fig. 6**. Numerical results in Example 2: the distributions of (a) analytical solution, (b) PIKFNN solution, and (c) relative error.

Table 11 presents the $L_2$ relative errors with different parameter settings, it can be found that the increasing training data will increase the predicting accuracy of the PIKFNN, and with the increasing wavenumber, the predicting accuracy will decrease, however, if the training data is further increased, the PIKFNN will give a better result. Fig. 6 plots the numerical results with the wavenumber $k=150$, the $tol=1E-4$, and 400 training data, we can find that the PIKFNN obtains a very accurate result in the whole computational domain. In general, the PIKFNN is feasible to simulate the wave propagation with high wavenumber.

**Example 3.** Laplace equation in an infinite domain.

Consider the following Laplace equation in an infinite domain as shown in Fig. 7,

$$\Delta u(\mathbf{x}) = 0, \quad \mathbf{x} = (x_1, x_2) \in \Omega, \tag{26}$$

$$u(\mathbf{x}) = \frac{x_1 + x_1}{x_1^2 + x_2^2}, \quad \mathbf{x} = (x_1, x_2) \in \partial\Omega. \tag{27}$$

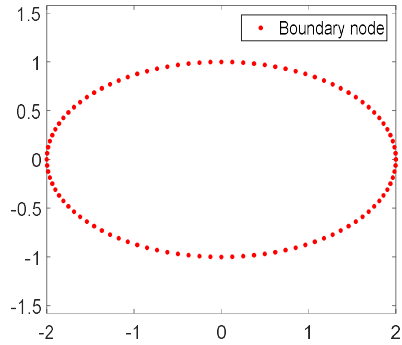

**Fig. 7.** Schematic diagram of the discretization of the infinite domain in example 3.

The schematic diagram of the discretization is depicted in Fig. 7. The number of boundary nodes is $N_B = 1600$, that is the training set contains $N = N_B = 1600$ samples, the source nodes are uniformly placed on a circle of radial 0.5 centered at the origin. The number of PIKFs was chosen to be $N = 1600$. The PIKF corresponding to the partial differential operator $\Delta$ is $\varphi^{L_0} = -\frac{1}{2\pi}\ln(r)$. To minimize the loss function, the L-M algorithm is employed in this example, and the $tol$ is set as 1E-10.

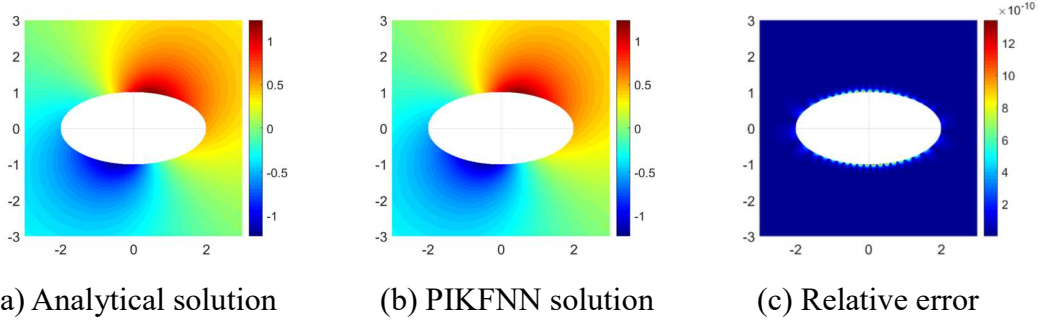

(a) Analytical solution     (b) PIKFNN solution     (c) Relative error

**Fig. 8**. Numerical results in Example 3: the distributions of (a) analytical solution, (b) PIKFNN solution, and (c) relative error.

Table 12 $L_2$ relative errors with different *tol* values in example 3.

| *tol* | 1E-4 | 1E-6 | 1E-8 | 1E-10 |
|---|---|---|---|---|
| $L_2$ relative error | 1.27E-7 | 4.13E-8 | 7.17E-9 | 7.94E-12 |
| CPU time (s) | 1.04 | 1.13 | 1.27 | 1.46 |

Fig. 8 shows the numerical results obtained by the proposed PIKFNN, it can be found that the PIKFNN solution is in good agreement with the analytical solution, and the relative error is at level 1E-10. Table 12 presents the $L_2$ relative errors with different *tol* values, we can find that the numerical error is increasing with the increasing *tol* value, correspondingly the CPU time will decrease. In general, the PIKFNN is feasible to solve boundary value problems in infinite computational domains without any additional artificial boundary conditions, due to the PIKFs naturally satisfy the boundary conditions at infinity.

**Example 4.** 3D nonhomogeneous modified Helmholtz equation.

Consider the following nonhomogeneous modified Helmholtz equation in a rabbit model as shown in Fig. 9,

$$(\Delta-1)u(\mathbf{x}) = 2e^{x_1+x_2+x_3}, \quad \mathbf{x}=(x_1, x_2) \in \Omega, \quad (28)$$

$$u(\mathbf{x}) = e^{x_1+x_2+x_3}, \quad \mathbf{x}=(x_1, x_2) \in \partial\Omega. \quad (29)$$

The schematic diagram of the discretization of the rabbit model is depicted in Fig. 9. The number of boundary nodes is $N_B = 1154$, the source nodes are uniformly placed on a sphere of radial 3 centered at the origin. The number of PIKFs was chosen to be

$2N$. The PIKFs corresponding to the partial differential operator $(\Delta-1)$ is $\varphi^{L_0}=\dfrac{e^{-r}}{4\pi r}$. The partial differential operator $(\Delta-3)$ is constructed to eliminate the nonhomogeneous term and the PIKFs corresponding to the partial differential operator $(\Delta-3)$ is $\varphi^{L_1}=\dfrac{e^{-\sqrt{3}r}}{4\pi r}$.

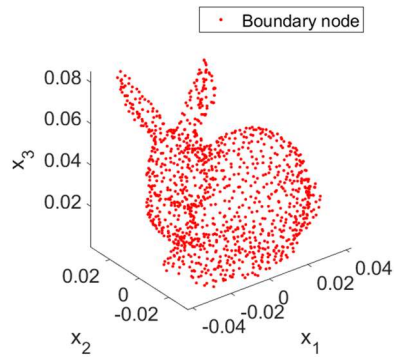

**Fig. 9.** Schematic diagram of the discretization of the rabbit model in example 4.

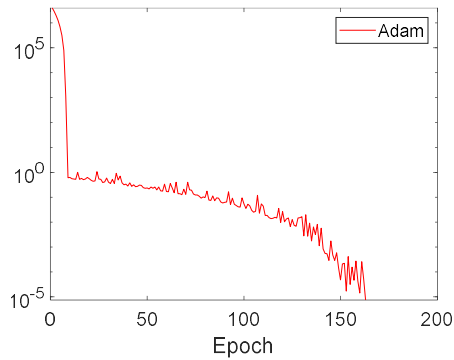

**Fig. 10.** Loss value varied with training epoch in example 4.

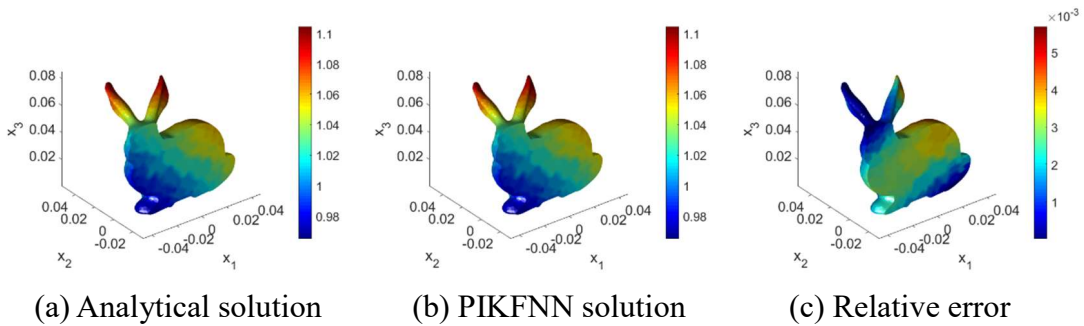

(a) Analytical solution      (b) PIKFNN solution      (c) Relative error

**Fig. 11**. Numerical results in Example 4: the distributions of (a) analytical solution, (b) PIKFNN solution, and (c) relative error.

Fig. 10 shows the loss value varied with the increasing training epoch, after 162 epochs, the loss value goes to the goal 1E-5, and Fig. 11 depicts the numerical results obtained by the PIKFNN, the relative error is at level 1.0E-3. This example demonstrates that the PIKFNN is feasible to solving nonhomogeneous PDEs in 3D complicated computational domains.

**Example 5.** 3D long-term heat conduction analysis under time-dependent thermal loading.

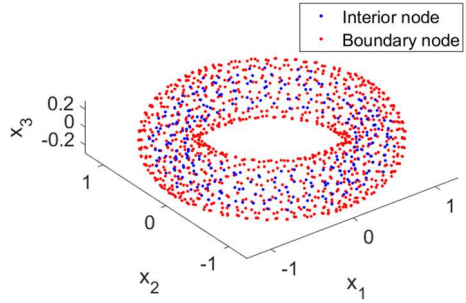

**Fig. 12.** Schematic diagram of the discretization of the torus model in example 5.

Consider the following transient heat conduction problem in a torus model as shown in Fig. 12,

$$\frac{\partial u(\mathbf{x},t)}{\partial t} - 0.001\Delta u(\mathbf{x},t) = f(\mathbf{x},t), \quad \mathbf{x} \in \Omega, \ 0 < t \leq T, \tag{30}$$

$$u(\mathbf{x},t) = \left(\sin(x_1) + \cos(x_2) + \sin(x_3)\right)e^{-0.003t}, \quad \mathbf{x} \in \partial\Omega, \ 0 < t \leq T, \tag{31}$$

$$u(\mathbf{x},0) = \left(\sin(x_1) + \cos(x_2) + \sin(x_3)\right), \quad \mathbf{x} \in \Omega. \tag{32}$$

In which $f(\mathbf{x},t) = -0.002\left(\sin(x_1) + \cos(x_2) + \sin(x_3)\right)e^{-0.003t}$ and the final time is $T = 100\ s$. The number of discrete nodes of the torus model is $1279$, which contains 998 boundary nodes and 281 interior nodes, where the initial condition number is $N_{Ini} = 1279$ and the boundary condition number is $N_{Bou} = 998 \times 5$ (uniformly distributed at the time instants $(20,\ 40,\ 60,\ 80,\ 100)s$). The spatial-temporal source node pairs have the same locations as the boundary and initial condition node pairs but have a delay time $dt = 200\ s$ in the temporal domain. The number of PIKFs was chosen to be $2(N_{Bou}+N_{Ini})=2N$. The PIKF corresponding to the spatial-temporal differential operator $\left(\dfrac{\partial}{\partial t}-0.001\Delta\right)$ is $\varphi^{\bar{L}_0}=\Theta(t-\tau)\dfrac{e^{r^2/0.004(t-\tau)}}{(0.004\pi(t-\tau))^{3/2}}$, and the spatial-temporal differential operator $\left(\dfrac{\partial}{\partial t}-0.002\Delta\right)$ is constructed to eliminate the thermal loading $f(\mathbf{x},t)$ and the corresponding PIKF is $\varphi^{\bar{L}_1}=\Theta(t-\tau)\dfrac{e^{R^2/0.008(t-\tau)}}{(0.008\pi(t-\tau))^{3/2}}$. To minimize the loss function, the L-M algorithm is employed in this example, the *tol* is set as 1E-10.

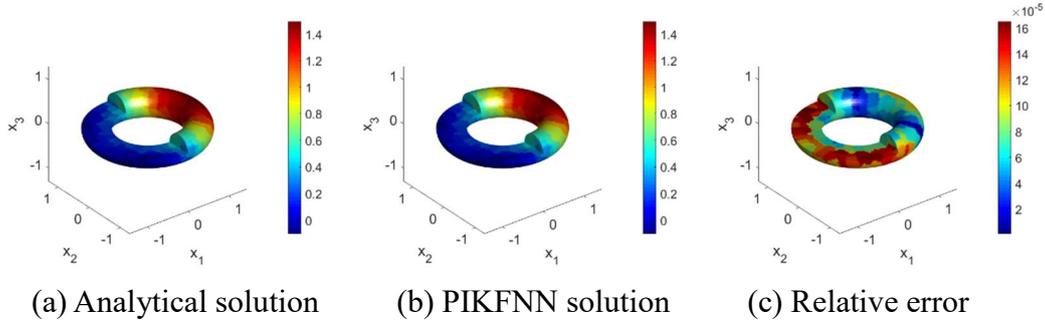

(a) Analytical solution    (b) PIKFNN solution    (c) Relative error

**Fig. 13**. Numerical results at $T=100\ s$ in Example 5: the distributions of (a) analytical solution, (b) PIKFNN solution, and (c) relative error.

Fig. 13 plots the numerical results at $T=100\ s$ obtained by the PIKFNN, it can be found that the PIKFNN solution is in good agreement with the analytical solution, where Fig. 13(c) plots the relative error distribution, the maximin relative error is 1.63E-4. This example verifies that the proposed PIKFNN is feasible and accurate to simulate the long-term transient heat conduction behavior. Moreover, it should be pointed out that the proposed PIKFNN also can be easily used to solve more general time-dependent problems.

**Example 6.** 2D spatial structural derivative diffusion model.

Consider the following spatial structural derivative diffusion model,

$$\dfrac{\partial u(\mathbf{x},t)}{\partial t}=k\left(\dfrac{d}{d_s x_1}\left(\dfrac{du(x_1,t)}{d_s x_1}\right)+\dfrac{d}{d_s x_2}\left(\dfrac{du(x_2,t)}{d_s x_2}\right)\right),\quad \mathbf{x}\in\Omega, 0<t\leq T, \qquad (33)$$

$$u(\mathbf{x}, t=0) = \cos\left(\frac{\pi f(x_1)}{2}\right) + \cos\left(\frac{\pi f(x_2)}{2}\right) + \sin\left(\frac{\pi f(x_1)}{2}\right) + \sin\left(\frac{\pi f(x_2)}{2}\right), \quad \mathbf{x} \in \Omega,$$

(34)

$$u(\mathbf{x}, t) = \left(\begin{array}{c} \cos\left(\frac{\pi f(x_1)}{2}\right) + \cos\left(\frac{\pi f(x_2)}{2}\right) \\ + \sin\left(\frac{\pi f(x_1)}{2}\right) + \sin\left(\frac{\pi f(x_2)}{2}\right) \end{array}\right) e^{\left(\frac{-\pi^2 kt}{4}\right)}, \quad \mathbf{x} \in \partial\Omega, 0 < t \leq T. \quad (35)$$

In which, $k$ denotes the diffusion coefficient, the final time is $T = 1\,s$, $f(x_i)$ are the spatial kernel functions, in this example the following four cases are investigated: $f(x_i) = x_i^{0.5}$, $f(x_i) = x_i$, $f(x_i) = e^{x_i}$, $f(x_i) = \ln x_i$. The 2D geometry and the (2+1)D computational domain $\Omega$ with discrete spatial-temporal nodes are depicted in Fig. 14. The number of the discrete spatial-temporal nodes of the (2+1)D computational domain $\Omega$ is $N = 5758$, and $N = 5758$ boundary and initial condition values are set as the targets, and their spatial-temporal coordinates are set as the labels. The spatial-temporal source node pairs are the same spatial-temporal nodes as the discrete spatial-temporal nodes but have a delay time $dt = 3\,s$ in the temporal domain. The number of PIKFs was chosen to be $N = 5758$. The PIKF corresponding to the spatial-temporal differential operator is $\varphi^{\bar{L}_0} = \Theta(t - \tau) \frac{e^{-(f(x_1) - f(x_2))^2 / (t - \tau)}}{4\pi(t - \tau)}$. To minimize the loss function, the L-M algorithm is employed in this example, the *tol* is set as 1E-10.

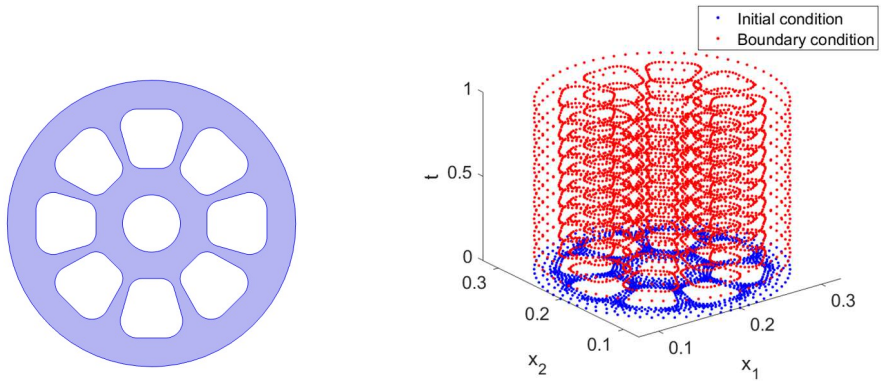

(a) Geometry　　　　　(b) Distributions of boundary and initial conditions

**Fig. 14**. Geometry and the schematic diagram of the training data in example 6: (a) the geometry, (b) the distribution of the training data, red node represents the boundary

condition and blue node represents the initial condition.

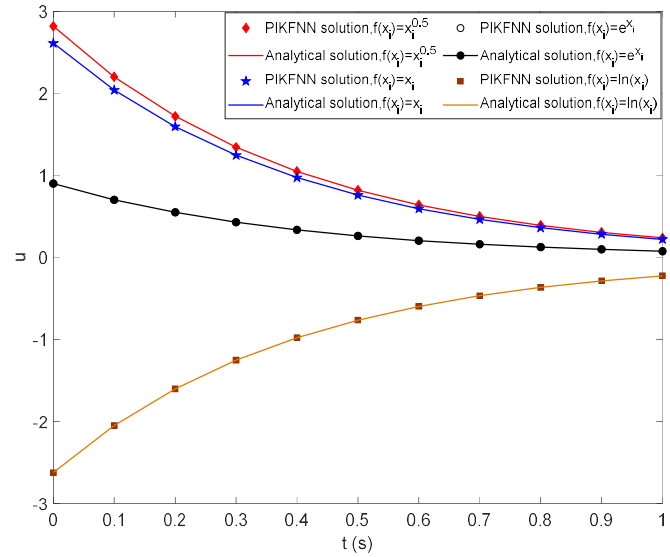

**Fig. 15**. Comparison of time evolutions of the PIKFNN solution at the point (0.31, 0.2) with different spatial kernels $f(x_i)$.

As can be seen from Fig. 15, the numerical solutions agree well with the analytical solutions under different spatial kernels, and the solution variations at the point (0.31, 0.2) show that the different kernels could lead to different variations trends of solution. Fig. 16 shows the numerical results obtained by the proposed PIKFNN at $T = 1\,s$, it can be found that the spatial kernel $f(x_i)$ has a great influence on the diffusion process. This also indicates that different anomalous diffusion phenomena could be described by the spatial structural derivative diffusion models with different spatial kernels.

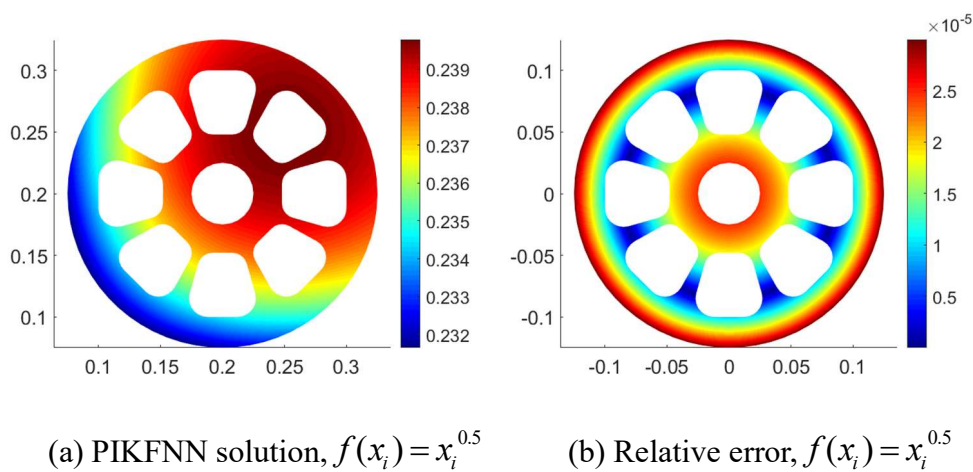

(a) PIKFNN solution, $f(x_i) = x_i^{0.5}$       (b) Relative error, $f(x_i) = x_i^{0.5}$

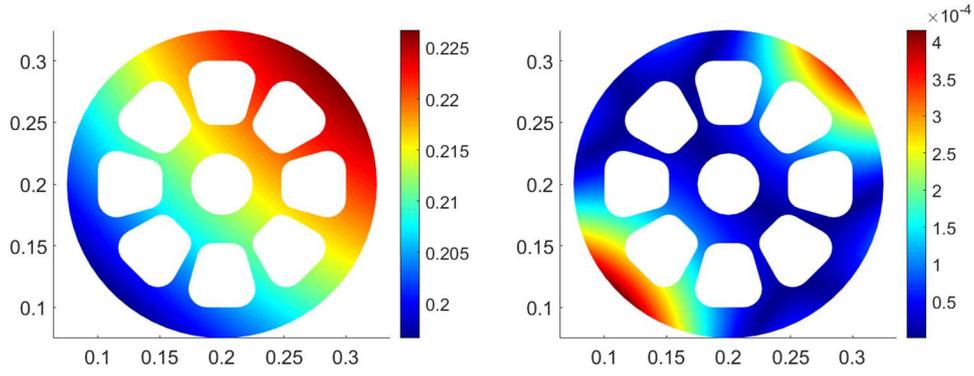

(c) PIKFNN solution, $f(x_i) = x_i$     (d) Relative error, $f(x_i) = x_i$

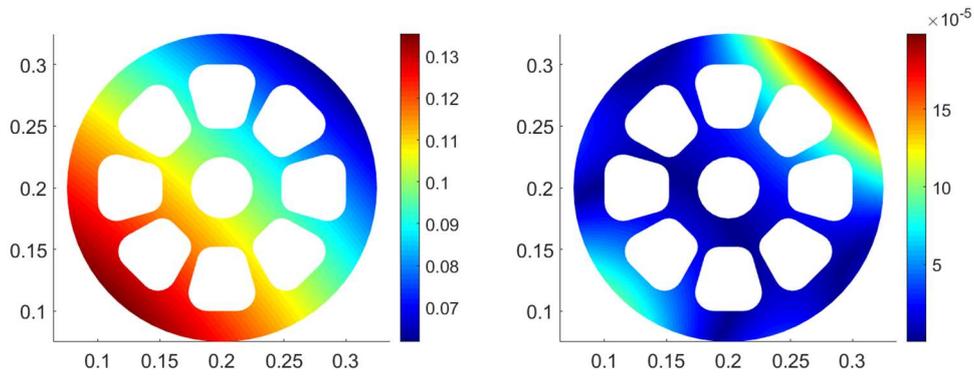

(e) PIKFNN solution, $f(x_i) = e^{x_i}$     (f) Relative error, $f(x_i) = e^{x_i}$

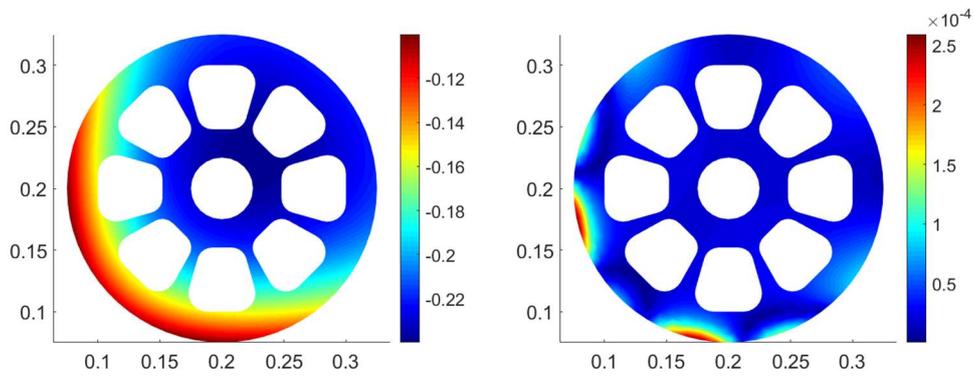

(g) PIKFNN solution, $f(x_i) = \ln x_i$     (h) Relative error, $f(x_i) = \ln x_i$

**Fig. 16**. Numerical results obtained by the proposed PIKFNN at $T = 1\ s$ in Example 6: PIKFNN solution distributions and relative error distributions with different spatial kernels, (a) $f(x_i) = x_i^{0.5}$, (b) $f(x_i) = x_i$, (c) $f(x_i) = e^{x_i}$, (d) $f(x_i) = \ln x_i$.

**Example 7.** High dimension Laplace equation.

In this case, we consider the following 4D Laplace equation,

$$\Delta u(\pmb{x}) = 0, \pmb{x} = (x_1, x_2, x_3, x_4) \in \Omega, \tag{36}$$

with Dirichlet boundary condition,

$$u(\pmb{x}) = x_1^2 + x_2^2 - x_3^2 - x_4^2, \pmb{x} = (x_1, x_2, x_3, x_4) \in \partial\Omega. \tag{37}$$

The computational domain $\Omega$ is a 4D hypersphere which is defined by

$$x_1^2 + x_2^2 + x_3^2 + x_4^2 \leq 1. \tag{38}$$

The number of the discrete nodes of the 4D computational domain $\Omega$ is $N = 1630$, and $N_B = 1630$ Dirichlet boundary condition values are set as the targets, and their coordinates are set as the labels. The source nodes are distributed at the boundary of the 4D hypersphere with radius 5 which is centered at the origin. The number of PIKFs was chosen to be $N = 1630$ The PIKF corresponding to the 4D Laplace operator is $\varphi^{L_0} = \frac{1}{8V_4}\frac{1}{r^2}$, where $V_4 = \frac{\sqrt{\pi^4}}{2}1^4$ is the volume of a unit 4D hypersphere. To minimize the loss function, the L-M algorithm is employed in this example, the *tol* is set as 1E-8.

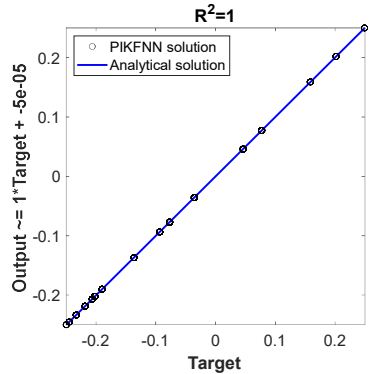

**Fig. 17**. Linear regression between PIKFNN solution and analytical solution.

To verify the accuracy of the PIKFNN for solving the high dimension Laplace equation, the PIKFNN solution at the surface of a 4D hypersphere with radius 0.5 is checked, which contains $Nt = 1630$ test points. Fig. 17 plots the linear regression between the PIKFNN solution and analytical solution at these test points, it can be found that the PIKFNN solution is in good agreement with the analytical solution, and the determination coefficient $R^2$ is equal to 1, which means the PIKFNN gives a very accurate prediction. Generally, this example demonstrated that the PIKFNN could solve High dimension PDE accurately.

**Example 8.** Potential-based inverse electromyography of pregnant sheep uterine contractions.

Based on the electro-quasi-static assumption (Wu, et al., 2019), the inductive and capacitance effects of the volume conductor of the pregnant sheep are ignored, thus the electrical potentials at each time instant satisfy the 3D Laplace equation. The electrical potentials on the sheep's body surface can be measured by some artificially placed electrodes, which are defined as the Dirichlet boundary conditions. And the electrical fluxes on the sheep's body surface are considered to be zero, due to the pregnant sheep's body is an insulating medium in the air, which are defined as the Neumann boundary conditions.

The numerical model of potential-based inverse electromyography can be considered as an inverse Cauchy problem with the over-specified boundary conditions on the sheep's body surface, via solving the above Cauchy problem, the electrical potentials on pregnant sheep's uterine surface at each time frame can be computed noninvasively.

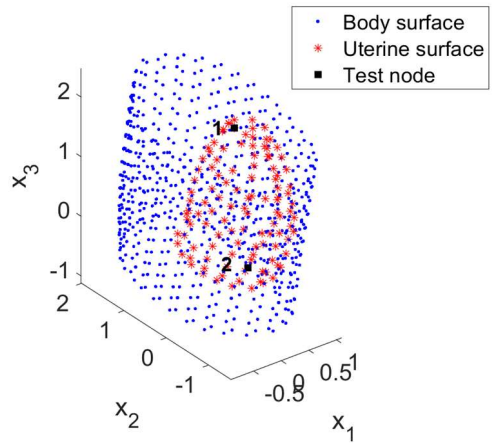

**Fig. 18**. Schematic diagram of the discretization of the sheep's body in example 8: Blue nodes represent the electrodes on the body surface, red nodes represent the discrete nodes on the uterine surface, and dark nodes represent the test nodes.

The number of the discrete nodes of the sheep's body is $N = 783$ as shown in Fig. 18, where 669 electrodes are placed on the body surface and the uterine surface was discretized into 114 interior nodes. Due to the over-specified boundary conditions on the sheep's body surface, the training set contains $2N$ samples, that is 669 Dirichlet boundary condition values and 669 Dirichlet boundary condition values are set as the targets, and their coordinates are set as the labels. The number of PIKFs was chosen to be $N = 783$. The PIKF corresponding to the 3D Laplace operator is $\varphi^{L_0} = \frac{1}{4\pi r}$. To minimize the loss function, the L-M algorithm is employed in this example, the *tol* is set as 1E-8.

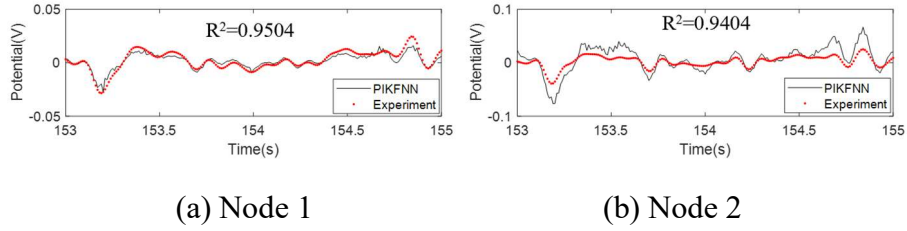

(a) Node 1  (b) Node 2

**Fig. 19**. Predicted electromyography at two electrodes on the uterine surface via the PIKFNN in comparison with the measured electromyography (Wu, et al., 2019) from 153 *s* to 155 *s*.

The predicted electromyography at two electrodes (node 1 (0.43, 0.13, 1.87), node 2 (−0.02, −0.77, 0.10)) on the uterine surface by the PIKFNN is compared to the measured electromyography from 153 *s* to 155 *s* with a time interval 0.01 *s*, which are displayed in Fig. 19. It can be found from Fig.19 that the predicted electromyography is in good agreement with the measured electromyography, the correlation coefficient $R^2$ between the predicted potentials and the measured potentials are pretty high. In general, the proposed PIKFNN could be considered a feasible approach in the application of potential-based inverse electromyography.

**Example 9.** Nonhomogeneous PDE with enhanced PIKFs.

Consider the following nonhomogeneous Helmholtz equation in a L-shaped domain as shown in Fig. 20,

$$(\Delta+1)u(\mathbf{x}) = x_1, \quad \mathbf{x} = (x_1, x_2) \in \Omega, \tag{39}$$

$$u(\mathbf{x}) = \sin(x_1) + \sin(x_2) + x_1, \quad \mathbf{x} = (x_1, x_2) \in \partial\Omega. \tag{40}$$

In the above examples, additional PIKF neurons centered at the boundary nodes are added for simulating the nonhomogeneous PDEs. In this example, we proposed an alternative approach for solving the nonhomogeneous PDEs, the enhanced PIKF neurons centered at the interior nodes are instead of the aforementioned PIKF neurons centered at the boundary nodes, and all the PIKFs are set as the same type, but the radial $r$ in the PIKF are modified to $\sqrt{r^2+s^2}$, where $s$ is a constant in this example. Also, due to the enhanced PIKF is no longer satisfies the governing equation, the loss function modified to

$$Loss_m(w) = \frac{1}{N_B}\sum_{i=1}^{N_B}\left(u_{PIKFNN} - \left(\sin(x_{1i}) + \sin(x_{2i}) + x_{1i}\right)\right)^2 + \frac{1}{N_I}\sum_{i=1}^{N_B}\left((\Delta+1)u_{PIKFNN} - x_{1i}\right)^2,$$

(41)

where the number of boundary nodes is $N_B = 62$ and interior nodes is $N_I = 209$. The number of PIKFs was chosen to be $N = N_B + N_I = 271$. The PIKFs is $\varphi^{L_0} = \frac{1}{2\pi}Y_0(\sqrt{r^2+s^2})$. To minimize the loss function $Loss_m$, the L-M algorithm is employed in this example, the *tol* is set as 1E−10.

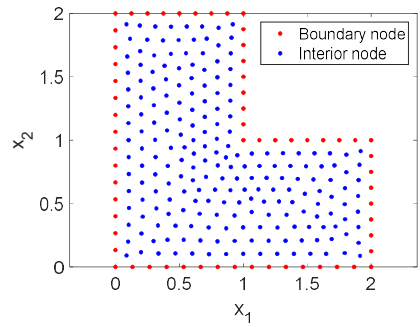

**Fig. 20**. Schematic diagram of the discretization of the L-shaped domain in example 9.

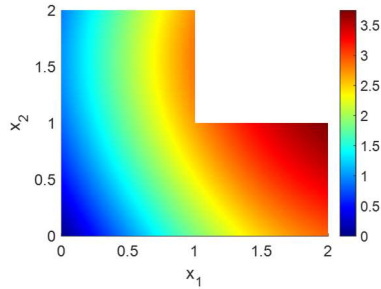

(a) Analytical solution, $s$=0.5

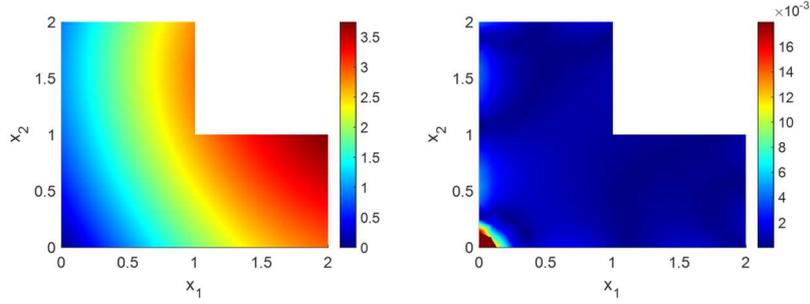

(b) PIKFNN solution, *s*=0.5   (c) Relative error, *s*=0.5

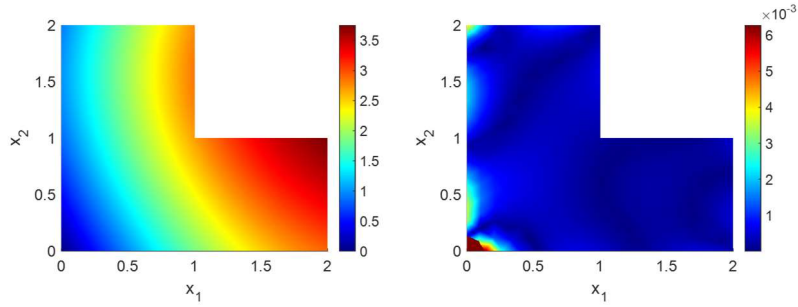

(b) PIKFNN solution, *s*=1   (c) Relative error, *s*=1

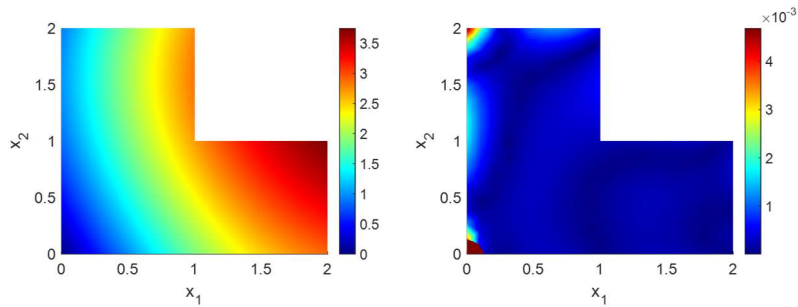

(b) PIKFNN solution, *s*=2   (c) Relative error, *s*=2

**Fig. 21**. Numerical results in Example 9: (a) the distributions of analytical solution; the distributions of (b) PIKFNN solution and (c) relative error with *s*=0.5; the distributions of (d) PIKFNN solution and (e) relative error with *s*=1; the distributions of (f) PIKFNN solution and (g) relative error with *s*=2.

Fig. 21 plots the numerical results of the nonhomogeneous Helmholtz equation (27)-(28) obtained by the PIKFNN, it can be found that the PIKFNN solutions with different constants are in good agreement with the analytical solution, and all the relative errors with different constant *s* are at level 1E−3. This example verifies that the proposed approach for solving nonhomogeneous PDEs is feasible and accurate.

**Example 10.** Elastic thin plate subjected to a uniform pressure.

In this example, we consider the following elastic thin plate under uniform pressure $p=1\ N/m^2$, as shown in Fig. 22, where the displacement components in both the $x_1$ and $x_2$ are constrained on the bottom boundary. And the material constants are taken to be the Poisson's ratio $v=0.3$ and the elastic modulus $\mu = 384615\ Pa$. The plane strain assumption is adopted in this problem since we assume the length of plate in $x_3$ direction is very large than other directions.

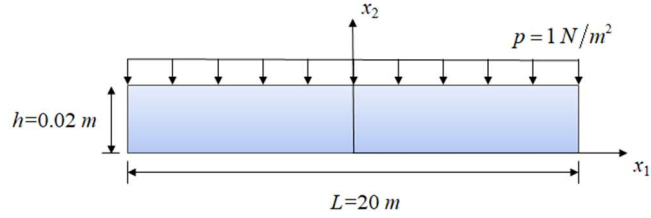

**Fig. 22.** Elastic thin plate subjected to uniform pressure.

Consider the following governing equation and analytical solution in this problem,

$$\sigma_{ij,j} = 0,\ i,j=1,2, \tag{42}$$

$$u_2(x_1,x_2) = -\frac{(1+v)(1-2v)}{2(1-v^2)\mu}x_2,\ \sigma_{11}=\frac{v}{v-1},\ \sigma_{22}=-1, \tag{43}$$

where $\sigma_{11}$ and $\sigma_{22}$ represent the normal stresses in the $x_1$ and $x_2$ directions, $\sigma_{12}$ and $\sigma_{21}$ denote the shear stress, $u_2$ stands for the displacement in $x_2$ direction. The number of boundary nodes is $N=108$, the training set contains $N=108$ samples, and the source nodes are uniformly placed on a sphere of radial 10 centered at the origin. The PIKFs corresponding to the displacement component and traction component of the governing equation are shown as the following:

$$\varphi_{disp}^{(\Re)} = \frac{1}{8\pi\mu(1-v)}\left\{(3-4v)\log\frac{1}{r}\delta_{l,k}+r_{,l}r_{,k}\right\},\ (l,k=1,2), \tag{44}$$

$$\varphi_{trac}^{(\Re)} = \frac{1}{4\pi(1-v)r}\{[(1-2v)\delta_{lk}+2r_{,l}r_{,k}]r_{,n}+(1-2v)(r_{,l}n_k-r_{,k}n_l)\}, \ (l,k=1,2), \quad (45)$$

where $r_{,n}=r_{,1}n_1+r_{,2}n_2$ expresses the derivative of $r$ in the direction of the outward normal, $n_i$ stands for the unit outward normal.

The numerical results obtained by the proposed PIKFNN are plotted in Fig. 23, we can observe that the displacement in $x_2$ direction computed by using PIKFNN is excellently consistent with the analytical solution for a small aspect ratio and with a relatively small number of training nodes.

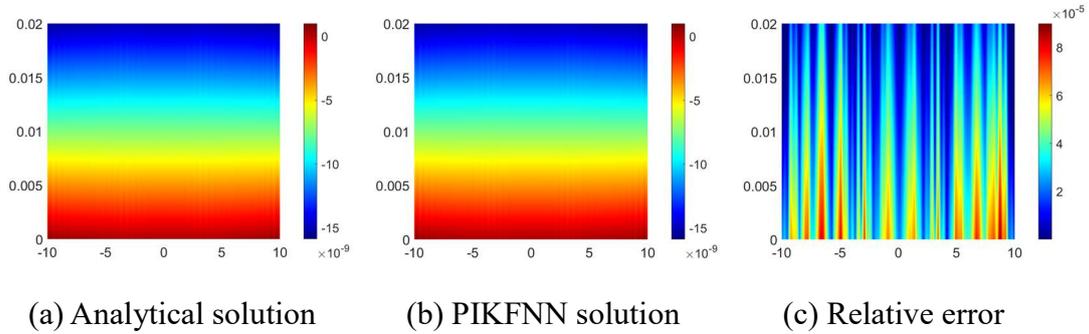

(a) Analytical solution  (b) PIKFNN solution  (c) Relative error

**Fig. 23**. Numerical results in Example: the distributions of (a) analytical solution, (b) PIKFNN solution, and (c) relative error.

Then, we fix the length of the elastic plate as $L = 20\ m$ while the thickness-to-length ratio experienced a decrease in the aspect ratio as the thickness reduces. Table 13 shows the relative errors of stress solutions at point (0, h/2), as the thickness-to-length ratio varies from $10^{-9}$ to $10^{-3}$. It can be observed that the predicted results calculated by using PIKFNN agree well with the analytical solutions, even though the aspect ratio h/L is very small. The results demonstrate the accuracy and validity of the proposed PIKFNN for thin elastic structures.

Table 13 Relative errors of stresses at the point (0, h/2).

| h/L | Stresses $\sigma_{11}$ | | Stresses $\sigma_{22}$ | |
|---|---|---|---|---|
| | Exact solution | Relative errors | Exact solution | Relative errors |
| 1E-03 | -0.4286 | $3.0697\times 10^{-10}$ | -1.0 | $1.4175\times 10^{-10}$ |
| 1E-04 | | $2.4524\times 10^{-8}$ | | $1.2092\times 10^{-9}$ |

| | | |
|---|---|---|
| 1E-05 | $7.3606\times10^{-9}$ | $1.0967\times10^{-9}$ |
| 1E-06 | $8.4756\times10^{-9}$ | $4.1915\times10^{-11}$ |
| 1E-07 | $2.9105\times10^{-9}$ | $1.2985\times10^{-10}$ |
| 1E-08 | $1.1601\times10^{-8}$ | $7.5609\times10^{-10}$ |
| 1E-09 | $2.8902\times10^{-8}$ | $1.2532\times10^{-9}$ |

## 5. Conclusions

In this paper, the use of physics-informed kernel function neural networks (PIKFNNs) for solving various linear PDEs has been explored, the physics-informed kernel functions which contain the prior physical information of considering PDEs were used as the activation functions of the PIKFNNs, the loss function which only includes the residual of the boundary/initial conditions was constructed. The numerical results show that the PIKFNNs converge to the decreasing loss value, and the PIKFNNs are feasible in solving long-time evolution problems, high wavenumber problems, inverse problems, and infinite domain problems without any additional technologies.

This study points out that introducing prior physical and mechanical information to activation functions in artificial neural networks will help them generalize well and overcome some challenges in the present physics-informed machine learning. The addition of prior physical and mechanical information can greatly reduce the neural network complexity, a two-layer PIKFNN can solve various linear PDEs. However, as the first step, there are still some challenges. Only the linear PDEs were solved in the present study, also, no theory guarantees that the present training procedure converges to a global minimum, which is an open problem to be solved.

**Declaration of competing interest**

The authors declare that they have no known competing financial interests or personal relationships that could have appeared to influence the work reported in this paper.

**Acknowledgments**

This work was funded by the National Natural Science Foundation of China (Grant

No. 12122205), the Six Talent Peaks Project in Jiangsu Province of China (Grant No. 2019-KTHY-009).